\documentclass[12pt,oneside, letterpaper]{article}
\usepackage{fancyhdr}  				
\usepackage{amssymb}					
\usepackage{amsmath}					
\usepackage{amsfonts}					
\usepackage{amsthm}					  
\usepackage{mathrsfs}					
\usepackage{lscape}						

\usepackage{amsthm, t1enc, amssymb, amsfonts, amsmath}
\usepackage{graphicx}
\usepackage[all]{xy}
\usepackage[dvips]{color}
\usepackage[active]{srcltx}
\usepackage{pstricks}
\usepackage[none]{hyphenat}

\usepackage[none]{hyphenat}  

\newcommand{\N}{\ensuremath{\mathbb{N}}}

\theoremstyle{plain}
\newtheorem{prop}{Proposition}

\newtheorem{theorem}[prop]{Theorem}

\newtheorem{definition}[prop]{Definition}

\begin{document}

\title{Intuitionistic Existential Graphs from a non traditional point of view.}

\author{Yuri A. Poveda  \and Steven Zuluaga  }

\maketitle



\renewcommand{\sectionmark}[1]{\markright{\thesection\ #1}}
\fancyhf{}
\fancyhead[LE,RO]{\scriptsize\thepage}
\fancyhead[LO]{\scriptsize\rightmark}
\fancyhead[RE]{\scriptsize\leftmark}
\renewcommand{\headrulewidth}{0.5pt}
\renewcommand{\footrulewidth}{0pt}
\addtolength{\headheight}{0.5pt}
\fancypagestyle{plain}{%
\fancyhead{}
\renewcommand{\headrulewidth}{0pt}
}

\sloppy   



{\setlength{\parindent}{20pt}



\textbf{Abstract.} 
In this article we develop a new version of the intuitionist existential graphs presented by Arnol Oostra \cite{Oostra}. The deductive rules presented in this article have the same meaning as those described in the work of Yuri Poveda \cite{yuri}, because the deductions according to the parity of the cuts are eliminated and are replaced by a finite set of recursive rules. This way, $ Alfa_I $ the existential graphs system for intuitional propositional logic follows the course of the deductive rules of the system $ Alfa_0 $ described by Poveda \cite{yuri}, and is equivalent to the intuitionistic propositional calculus. 

In this representation the $ Alfa_0 $ system is improved, there are a series of deductive rules of second degree incorporated that previously had not been considered and that allow a better management of deductions and finally from the ideas proposed by Van Dalen \cite{ Van}, a mixture is incorporated in the deduction techniques, the natural deductions of the Gentzen system are combined with new system rules $ Alfa_0 $ and $ Alfa_I $.

The symbols proposed for the  $Alfa_I$ representation relate open, closed and quasi-open sets of the usual topology of the plot with the intuitional propositional logic, usefull for approaching new problems in the representation of this logic from a more geometrical perspective.\\
 
\textbf{Keywords:} Propositional calculus, intuitionism, existential graphs, deductive rules.

}

\section{Introduction.}
The following assignment is an addition to the studies made on Peirce's existential graphs in Colombia, motivated by professor Yuri Alexander Poveda's suggestion of finding a system of existential graphs equivalent to the intuitionistic propositional calculus. Problem to which a solution is given in this article. 

Peirce's existential graphs are deductive systems that formalize the classic logic of propositions ($alfa$ system), of predicates ($beta$ system) and the modal ($gama$ system). Likewise, these are innovative systems that contribute a totally new vision of the logic principles, which facilitate learning, management, discovery, and truth deductions, besides offering other advantages. Unfortunately, these haven't had a great reception by the logical community; however, Roberts and Zeman have studied Peirce's graphs in their doctorate thesis (where they precisely study the $alfa$, $beta$ and $gama$ systems); later on, Burch, Brady and Trimble; and finally in Colombia, Oostra \cite{Oostra}, Zalamea \cite{Zalamea} and Poveda \cite{yuri}.

The biggest interest of this article is first, to present the intuitionistic system of existential graphs ($GEI$) obtained from the graphic systems $ALFAo$, proposed by Yuri A. Poveda \cite{yuri} and the natural deductive system for intuitionism presented by Van Dalen \cite{Van}. Secondly, to show that it is equivalent to the Intuitionistic Propositional Calculus. 

In the first section the $ALFAo$ system is presented with some modifications, where the introduction of a new rule that marks the difference between $alfa$ and $ALFAo$, and a new notation that facilitates making deductions in $ALFAo$, stand out. Likewise, the calculus of rules is formalized through four definitions that are introduced at the beginning of the chapter, the symmetries of the $ALFAo$ rules are studied, and the relations of these with the Modus Ponendo Ponens, and with the insertion and elimination of the double cut. In the second section the following things are done: a presentation of the natural deduction system for intuition, extracted from the Van Dalen \cite{Van} article; the definition of two new graphs used to design the implication and the disjunction of a new system, as well as a function used to translate graphs into formulas; the presentation of a system of $GEI$ ($ALFA_I$), who's deductive rules are obtained from the natural deduction system presented before; the change of some of $ALFA_I$ basic rules for other ones (theorems of the same); and the presentation of the $ALFA_{Io}$ graphs system that is equivalent to $ALFA_I$. Finally, a new not-intuitionistic rule is added to the $ALFA_{Io}$ system and the equivalence of this with $ALFAo$ is studied. In the third chpater: We approach exegetically the article of Arnold Oostra in which he exposes a $GEI$ system  \cite{Oostra} that, different to the $ALFA_{Io}$, harmonizes with $alfa$ conserving the same structure to enunciate the graphical transformation rules and the notation to present the demonstrations, and some differences between $ALFA_{Io}$ and the system proposed by Arnold Oostra are shown. 

\section{$ALFAo$ deductive system}

$ALFAo$ is a deductive system of graphs equivalent to the Classic Propositional Calculus, defined from the existential graphs system, $alfa$ by Charles Sander Peirce. Next we will present $ALFAo$ with some modifications that attend to a simpler presentation of the same, to see a more detailed presentation of this system, refer to \cite{Oostra}. 

\begin{flushleft}\textbf{Primitive systems:}\end{flushleft}

\begin{flushleft}
The rectangle $\xy(0,0)*=<0.5cm,0.75cm>\frm{-}\endxy$ , the closed curve  $\xy*\xycircle<0.25cm,0.15cm>{}\endxy$,  the letters $p_i$,$q_i$,$r_i$;$A_i$,$B_i$,$C_i$, $i\in \N$ and the letter $\lambda$.
\end{flushleft}

The following definition formalizes in a general way the types of rules that may take place in the system. In other systems like in \cite{Oostra} and \cite{Zalamea}, the definition of what a deductive rule limited to first degree rules is, appears indirectly. 

\begin{flushleft}\textbf{Deductive rules (rds):}\end{flushleft}

\begin{enumerate}
	\item Given $G_1$ y $G_2$ wdgs then  $G_1$$\vdash$$G_2$ is a first degree rd.
  \item Given $R_1$,$R_2$,...,$R_j$,$R$ first degree rds then $\xy(0,0)*=<2cm,0.01cm>\frm{-};(0,3)*+{R_1,R_2,...,R_j};(0,-3)*+{R}\endxy$ is a second degree rd. 
\end{enumerate}

\begin{flushleft}\textbf{System's axiom:}\end{flushleft}

\xy*=<0.5cm,0.75cm>\frm{-}\endxy  

\begin{flushleft}\textbf{System's deductive rules:}\end{flushleft}

The following rules are presented according to the deductive rule definition given previously, either of insertion or elimination. This classification appears indirectly in \cite{yuri}.

\begin{flushleft}\emph{A) First degree deductive rules:}\end{flushleft}

\begin{flushleft}Of insertion: A rule is of insertion if to pass from a graphic to another letters or closed curves called cuts are drawn.
\end{flushleft}

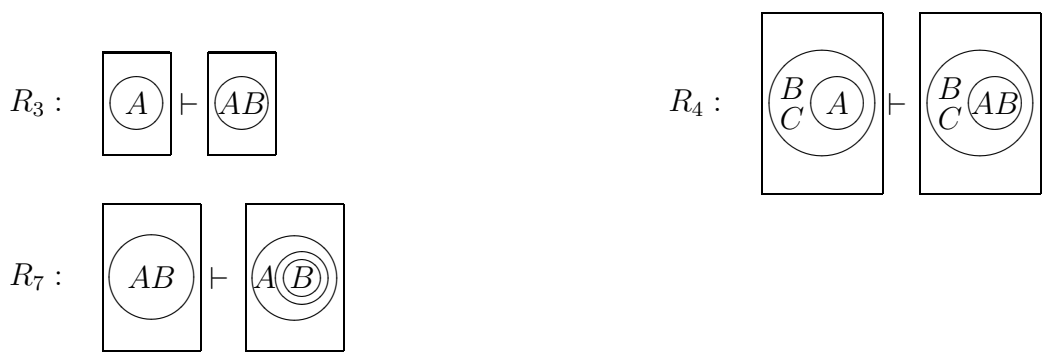
\begin{flalign*}
& R_3: \quad \xy(0,0)*=<0.9cm,1.35cm>\frm{-}*+{A}*\cir<10pt>{};(14,0)*=<0.9cm,1.35cm>\frm{-}*+{AB}*\cir<10pt>{};(7,0)*+{\vdash}\endxy
& R_4: \quad \xy(0,0)*=<1.6cm,2.4cm>\frm{-}*\cir<20pt>{};(2,0)*+{A}*\cir<10pt>{};(-4,2)*+{B};(-4,-2)*+{C};(21,0)*=<1.6cm,2.4cm>\frm{-}*\cir<20pt>{};(23,0)*+{AB}*\cir<10pt>{};(17,2)*+{B};(17,-2)*+{C};(10,0)*+{\vdash}\endxy \\
& R_7: \quad
\xy(0,0)*=<1.3cm,1.95cm>\frm{-}*+{AB}*\cir<16pt>{};(19,0)*=<1.3cm,1.95cm>\frm{-}*\cir<16pt>{};(15,0)*+{A};(20,0)*+{B}*\cir<7pt>{}*\cir<10pt>{};(9,0)*+{\vdash}\endxy
\end{flalign*}

\begin{flushleft}Of elimination: A rule is of elimination if to pass from a graphic to another letters or cuts are eliminated.
\end{flushleft}

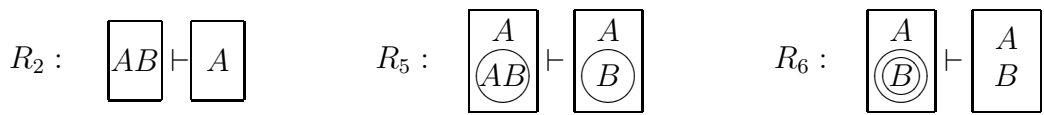
\begin{flalign*} & R_2: \quad \xy(0,0)*=<0.7cm,1.05cm>\frm{-}*+{AB};(11,0)*=<0.7cm,1.05cm>\frm{-}*+{A};(6,0)*+{\vdash}\endxy & R_5: \quad \xy(0,0)*=<0.9cm,1.35cm>\frm{-};(0,-2)*+{AB}*\cir<11pt>{};(14,0)*=<0.9cm,1.35cm>\frm{-};(14,-2)*+{B}*\cir<11pt>{};(0,4)*+{A};(14,4)*+{A};(7,0)*+{\vdash}\endxy & & R_6: \quad \xy(0,0)*=<0.9cm,1.35cm>\frm{-};(0,-2)*+{B}*\cir<7pt>{}*\cir<10pt>{};(0,4)*+{A};(14,0)*=<0.9cm,1.35cm>\frm{-};(14,-2)*+{B};(14,3)*+{A};(7,0)*+{\vdash}\endxy
\end{flalign*}

\begin{flushleft}\emph{B) Second degree deductive rule:}\end{flushleft}

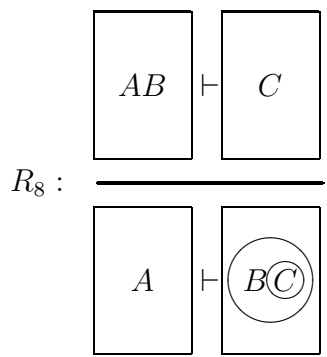
\begin{flalign*}&
\xy(-14,-13)*+{R_8:};(0,0)*=<1.3cm,1.95cm>\frm{-}*+{AB};(17,0)*=<1.3cm,1.95cm>\frm{-}*+{C};(9,0)*+{\vdash};(9,-13)*=<3cm,0.01cm>\frm{-};(0,-26)*=<1.3cm,1.95cm>\frm{-}*+{A};(17,-26)*=<1.3cm,1.95cm>\frm{-}*\cir<16pt>{};(15,-26)*+{B};(19,-26)*+{C}*\cir<7pt>{};(9,-26)*+{\vdash}\endxy \end{flalign*}

\subsection{Deduction in $ALFAo$}

\begin{definition}
$G\underset{ALFAo}{\vdash}G^{'}$ if and only if  exists $G_i\vdash G_{i+1}$ rds of $ALFAo$ with $1\leq i\leq n$ where $G=G_1$ y $G^{'}=G_{n+1}$ ($G^{'}$ is formal theorem if $G$ is the axiom or a formal theorem).
\end{definition}

\begin{definition}
Given the link $\{G_i \underset{ALFAo}{\vdash}G_{i+1}\}_{1\leq i\leq n}$ of rds then\\
$\xy(0,0)*=<3cm,0.01cm>\frm{-};(0,3)*+{\{G_i \underset{ALFAo}{\vdash}G_{i+1}\}_{1\leq i\leq n}};(0,-3)*+{G_1\vdash G_{n+1}}\endxy$ is a second degree rd of $ALFAo$.
\end{definition}

\begin{definition}
Given $\xy(0,0)*=<1cm,0.01cm>\frm{-};(0,3)*+{\{R_i\}_{1\leq i\leq n}};(0,-3)*+{R}\endxy$  and $\{R_i\}_{1\leq i\leq n}$ rds of $ALFAo$ then $R$ is a $ALFAo$ rd.
\end{definition}

Note: It can happen that  $\xy(0,0)*=<1cm,0.01cm>\frm{-};(0,3)*+{\{R_i\}_{1\leq i\leq n}};(0,-3)*+{R}\endxy$ be a $ALFAo$ rd and that the rds $\{R_i\}_{1\leq i\leq n}$ and $R$ not be of $ALFAo$. [see $R_8$].

\begin{definition}\label{rd2}
Given $\{R_i^{'}\}_{1\leq i\leq n}\notin ALFAo$ y $\{R_j\}_{1\leq j\leq m}\in ALFAo$ we have that
\begin{center}
$\xy(0,0)*=<1.5cm,0.01cm>\frm{-};(0,3)*+{\{R_i^{'}\}_{1\leq i\leq n}};(0,-3)*+{R}\endxy \Leftrightarrow \xy(0,0)*=<3cm,0.01cm>\frm{-};(0,3)*+{\{R_i^{'}\}_{1\leq i\leq n},\{R_j\}_{1\leq j\leq m}};(0,-3)*+{R}\endxy$
\end{center}
\end{definition}

The following rule, that translated to the Hilbert type systems correspond with the rule $\alpha\rightarrow\beta,\alpha\rightarrow\gamma \vdash \alpha\rightarrow\beta\wedge\gamma$, is used in the other existential graphic systems \cite{Oostra}, \cite{Zalamea}, \cite{Van} in an intuituve and informal way.
However here we introduce a rule of the system that isn't deductible from this one.
\vspace{10pt}

\xy(-10,-8)*+{R_0:};(12,0)*=<0.7cm,1.05cm>\frm{-}*+{B};(0,0)*=<0.7cm,1.05cm>\frm{-}*+{A};(6,0)*+{\vdash};(18,0)*+{,};(24,0)*=<0.7cm,1.05cm>\frm{-}*+{A};(36,0)*=<0.7cm,1.05cm>\frm{-}*+{C};(30,0)*+{\vdash};(18,-8)*=<4.5cm,0.01cm>\frm{-};(24,-16)*=<0.7cm,1.05cm>\frm{-}*+{A};(36,-16)*=<0.7cm,1.05cm>\frm{-}*+{BC};(30,-16)*+{\vdash}\endxy \vspace{10pt}

Introducing this rule in $ALFAo$ as rd of the system makes it possible to, with great ease, make the deduction of the rule $\alpha,\beta \vdash \alpha \wedge \beta$ known in the Hilbert type systems as insertion of $\wedge$.

\begin{theorem}\end{theorem}

\begin{center}
\xy(-10,-10)*+{R_0};(11,0)*=<0.7cm,1.05cm>\frm{-}*+{A};(0,0)*=<0.7cm,1.05cm>\frm{-};(6,0)*+{\vdash};(32,0)*=<0.7cm,1.05cm>\frm{-}*+{B};(21,0)*=<0.7cm,1.05cm>\frm{-};(27,0)*+{\vdash};(32,-20)*=<0.7cm,1.05cm>\frm{-}*+{AB};(21,-20)*=<0.7cm,1.05cm>\frm{-};(27,-20)*+{\vdash};(16,-10)*=<4cm,0.01cm>\frm{-};(16,0)*+{,}\endxy
\end{center}

The rule $R_1$ is deduced from $ALFAo$ therefore is a theorem of the system and can be suppressed as it's basic rule.

\begin{theorem}
$R_1:\xy(0,0)*=<0.7cm,1.05cm>\frm{-}*+{A};(11,0)*=<0.7cm,1.05cm>\frm{-}*+{AA};(6,0)*+{\vdash}\endxy \in ALFAo$
\end{theorem}

\begin{flalign*}&
\xy(12,0)*=<0.7cm,1.05cm>\frm{-}*+{A};(0,0)*=<0.7cm,1.05cm>\frm{-}*+{A};(6,0)*+{\underset{R_2}{\vdash}};(18,0)*+{,};(24,0)*=<0.7cm,1.05cm>\frm{-}*+{A};(36,0)*=<0.7cm,1.05cm>\frm{-}*+{A};(30,0)*+{\underset{R_2}{\vdash}};(18,-8)*=<4.5cm,0.01cm>\frm{-};(-8,-8)*+{R_0};(24,-16)*=<0.7cm,1.05cm>\frm{-}*+{A};(36,-16)*=<0.7cm,1.05cm>\frm{-}*+{AA};(30,-16)*+{\vdash}\endxy
\end{flalign*}

\begin{flushleft}
Note that by the definition \ref{rd2} $R_1 \in ALFAo$.
\end{flushleft}
Each of the basic rules of the $ALFAo$ system (of $R_1$ to $R_8$), seem to be necessary to form an equivalent system to $CPC$, this is because originally each and every one of these were defined based on the rules of Peirce's $alfa$ system. Nevertheless, the rules $R_7$ and $R_4$ can be supressed as the system's basic deductive rules, since they are theorems of $ALFAo$.

\begin{theorem}
$R_7$ is deductible by $R_5$,$R_2$ and $R_8$.
\end{theorem}

\begin{flalign*}&
\xy(22,15)*=<1.6cm,2.4cm>\frm{-};(22,10)*+{AB}*\cir<12pt>{};(22,20)*+{A};(43,15)*=<1.6cm,2.4cm>\frm{-};(43,10)*+{B}*\cir<12pt>{};(43,20)*+{A};(33,15)*+{\underset{R_5}{\vdash}};(64,15)*=<1.6cm,2.4cm>\frm{-};(64,10)*\cir<12pt>{}*+{B};(54,15)*+{\underset{R_2}{\vdash}};;(45,0)*=<5.5cm,0.01cm>\frm{-};(10,0)*+{R_8};(43,-15)*=<1.6cm,2.4cm>\frm{-}*+{AB}*\cir<20pt>{};(64,-15)*=<1.6cm,2.4cm>\frm{-}*\cir<20pt>{};(67,-15)*+{B}*\cir<7pt>{}*\cir<10pt>{};(61,-15)*+{A};(54,-15)*+{\vdash}\endxy
\end{flalign*}

\begin{theorem}
$R_4$ is deductible by $R_5$,$R_2$,$R_6$ and $R_8$.
\end{theorem}

\begin{theorem}
$R_8^-:\xy(0,13)*=<1.3cm,1.95cm>\frm{-}*+{A};(20,13)*=<1.3cm,1.95cm>\frm{-}*\cir<16pt>{};(16.5,13)*+{B};(21,13)*+{C}*\cir<7.5pt>{};(11,13)*+{\vdash};(11,0)*=<3.5cm,0.01cm>\frm{-};(0,-13)*=<1.3cm,1.95cm>\frm{-}*+{AB};(20,-13)*=<1.3cm,1.95cm>\frm{-}*+{C};(11,-13)*+{\vdash}\endxy \in ALFAo$
\end{theorem}

\begin{flalign*}&
\xy(-20,0)*=<1.3cm,1.95cm>\frm{-}*+{AB};(0,0)*=<1.3cm,1.95cm>\frm{-}*+{A};(20,0)*=<1.3cm,1.95cm>\frm{-}*\cir<16pt>{};(16.5,0)*+{B};(21,0)*+{C}*\cir<7.5pt>{};(11,0)*+{\underset{p}{\vdash}};(-11,0)*+{\underset{R_2}{\vdash}};(30,0)*+{,};(40,0)*=<1.3cm,1.95cm>\frm{-}*+{AB};(60,0)*=<1.3cm,1.95cm>\frm{-}*+{B};(51,0)*+{\underset{R_2}{\vdash}};(20,-13)*=<9.5cm,0.01cm>\frm{-};(-30,-13)*+{R_0};(40,-26)*=<1.3cm,1.95cm>\frm{-}*+{AB};(60,-26)*=<1.3cm,1.95cm>\frm{-}*\cir<16pt>{};(56.5,-26)*+{B};(61,-26)*+{C}*\cir<7.5pt>{};(57,-19)*+{B};(80,-26)*=<1.3cm,1.95cm>\frm{-}*+{C};(51,-26)*+{\vdash};(70,-26)*+{\underset{MP}{\vdash}}\endxy
\end{flalign*}

It is immediately noted that the strongest rules in $ALFAo$ are $R_8$ and $R_0$, without being coincidence that they are of second degree. Similarly, it is awaited that $R_8^-$ would be equally strong to it's inverse. This results to be true, which will be shown in the deductions of the inverse rules that remain.

\begin{prop}
$R_1^-,R_4^-,R_5^-,R_6^- y R_7^- \in ALAFAo$
\end{prop}

\subsection{$MP$ and the elimination of the double cut}

The Modus Ponendo Ponens is a theorem of $ALFAo$ 

\begin{theorem}The Modus Ponendo Ponens $MP:
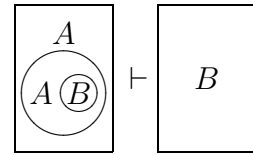
\xy(0,0)*=<1.3cm,1.95cm>\frm{-};(2,-2)*+{B}*\cir<7pt>{};(-3,-2)*+{A};(0,6)*+{A};(0,-2)*\cir<16pt>{};(10,0)*+{\vdash};(19,0)*=<1.3cm,1.95cm>\frm{-}*+{B}\endxy \in ALFAo$\end{theorem}

\begin{flalign*}&
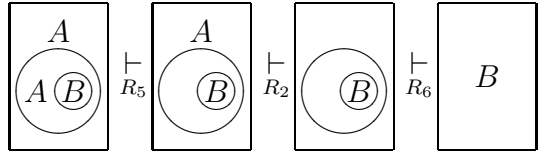
\xy(0,0)*=<1.3cm,1.95cm>\frm{-};(2,-2)*+{B}*\cir<7pt>{};(-3,-2)*+{A};(0,6)*+{A};(0,-2)*\cir<16pt>{};(10,0)*+{\underset{R_5}{\vdash}};(19,0)*=<1.3cm,1.95cm>\frm{-};(21,-2)*+{B}*\cir<7pt>{};(19,6)*+{A};(19,-2)*\cir<16pt>{};(29,0)*+{\underset{R_2}{\vdash}};(38,0)*=<1.3cm,1.95cm>\frm{-};(40,-2)*+{B}*\cir<7pt>{};(38,-2)*\cir<16pt>{};(48,0)*+{\underset{R_6}{\vdash}};(57,0)*=<1.3cm,1.95cm>\frm{-}*+{B}\endxy\end{flalign*}

The result presented below suggests a question: Can the Modus ponendo ponens, which is an intuitionist rule, deduce a non-intuitionist rule? The answer is clear, no. However, in $ALFAo$ the Modus Ponendo Ponens with the help of other rules, deduce the elimination of the double cut. For this reason, it is not possible to obtain an intuitionistic system of graphs from $ALFAo$, unless the vertex that unites these two rules is undone.

\begin{theorem}$MP$ deduces $R_6$.\end{theorem}

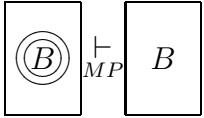
\begin{flushleft}
\xy(0,0)*=<1cm,1.5cm>\frm{-}*+{B}*\cir<7pt>{}*\cir<10pt>{};(8,0)*+{\underset{MP}{\vdash}};(16,0)*=<1cm,1.5cm>\frm{-}*+{B}\endxy
\end{flushleft}

\begin{flushleft}It is easy to see that with $R_2$,$R_0$ and $MP$ we can deduce $R_6.$\end{flushleft}

From the above, it can be stated that in $ALFAo$ the implication (and disjunction) is not independent of the negation and conjunction, like it does happen in the intuitionist propositional calculus ($IPC$). Relation that must be avoided to define the intuitionist system of graphs.

\section{Intuitionist Existential Graphs and the $IPC$}

The first presentation of a formal system for the intuitionist logic was published by Arend Heyting in 1930 and its formalization was presented in Hilbert style (two rules of inference and a large number of axioms). Four years later, Gerhard Gentzen announced two different alternatives to formalize it: the sequencing calculus and the natural deduction, the last characterized both by having rules of insertion And elimination for each connective, as by having an abbreviated notation for the deductions. Subsequently, Arnold Oostra presented an existential graphical version, equivalent to the $IPC$, symmetric in the original sense (the one adopted in the presentations for the $alpha$ system).\\

This chapter shows a $GEI$ system constructed from the natural deduction system for intuitionism presented by Dirk Van Dalen \cite{Van}, which preserves some of the characteristics that differentiate the $ALFAo$ system from the $alpha$ system.\\

The path chosen to find this system was totally different from that followed by Arnold Oostra, who was based primarily on the $alpha$ system and the Peircean legacy. For example, two new graphs were introduced, one to represent the implication and the other to represent the disjunction; However, with regard to denial and conjoint, $ALFAo$ graphs were used; In the case of the set of the basic deductive rules, we copied those of the natural deduction system for the aforementioned intuitionism, through a function that allows to make translations of formulas to graphs; And, finally, the graphical system obtained was refined.\\

\subsection{Intuitionist Existential Graphs}

Drawings that represent the implication and the disjunction: \\

$\xy(0,0)*\xycircle<20pt,18pt>{};(3,0)*\xycircle<7pt>{.}\endxy$ for the implication $\xy*\xycircle<20pt,18pt>{};(3,0)*\xycircle<7pt>{.}*\cir<7pt>{d^u};(-3,0)*\xycircle<7pt>{.}*\cir<7pt>{d^u}\endxy$ for the disjunction.

Note: the dotted curve was used by Peirce for the modal logic as a representation of the possible; However, here this drawing is taken with a totally different connotation that will be made known later.

\subsubsection{Natural Deduction for the intuitionism \cite{Van}.}

\begin{flushleft}
\textbf{Primitive symbols:}
\end{flushleft}

\begin{flushleft}
 The conjunction $\wedge$, the disjunction $\vee$, the implication $\rightarrow$, the parenthesis $($,$)$, the constants: false $\bot$ and true $\top$, and the letters $p_i$,$q_i$,$r_i$;$A_i$,$B_i$,$C_i$, $i\in \N$.
\end{flushleft}

\begin{definition}
well formed formulas ($fbf$):
\begin{enumerate}
	\item An atomic proposition $p$ is a $fbf$.
	\item The constants are $fbf$.
	\item If $A$ and $B$ are $fbf$s, then $A\wedge B$, $A\vee B$, $A\rightarrow B$ are $fbf$s.
	\item The formulas constructed according to 1. 2. y 3. are also $fbf$s.
\end{enumerate}
\end{definition}

\begin{definition}
$\neg$A=A$\rightarrow$$\bot$
\end{definition}

\textbf{System's deductive rules}

\emph{A) First degree deductive rules}\\

Of insertion

\begin{align*} & \wedge_i: \quad \xy(0,0)*=<2cm,0.01cm>\frm{-};(-4,3)*+{A};(4,3)*+{B};(0,-3)*+{A\wedge B}\endxy & \vee_i: \quad \xy(0,0)*=<1cm,0.01cm>\frm{-};(0,3)*+{A};(0,-3)*+{A\vee B}\endxy \quad \xy(0,0)*=<1cm,0.01cm>\frm{-};(0,3)*+{B};(0,-3)*+{A\vee B}\endxy
\end{align*}

Of elimination

\begin{align*} & \bot_{e}: \quad \xy(0,0)*=<1cm,0.01cm>\frm{-};(0,3)*+{\bot};(0,-3)*+{A}\endxy & \wedge_e: \quad \xy(0,0)*=<1cm,0.01cm>\frm{-};(0,3)*+{A\wedge B};(0,-3)*+{A}\endxy \quad \xy(0,0)*=<1cm,0.01cm>\frm{-};(0,3)*+{A\wedge B};(0,-3)*+{B}\endxy  & & \rightarrow_e: \quad \xy(0,0)*=<2cm,0.01cm>\frm{-};(-4,3)*+{A\rightarrow  B};(5,3)*+{A};(0,-3)*+{B}\endxy
\end{align*}

\begin{flushleft}
\emph{B) Second degree deductive rules}
\end{flushleft}

\begin{align*} & \rightarrow_i: \quad \xy(0,0)*=<1cm,0.01cm>\frm{-};(0,3)*+{B};(0,-3)*+{A\rightarrow B};(0,12)*+{A};(0,6)*=<0.5cm,0.01cm>\frm{-}\endxy[A] & \vee_e: \quad \xy(0,0)*=<3cm,0.01cm>\frm{-};(-10,3)*+{A\vee B};(0,3)*+{C};(0,12)*+{A};(0,7)*=<0.5cm,0.01cm>\frm{-};(10,3)*+{C};(10,12)*+{B};(10,7)*=<0.5cm,0.01cm>\frm{-};(0,-3)*+{C}\endxy[A][B]
\end{align*}

The formulas between brakets represent canceled premises in the new deduction.\\

The rule $\rightarrow_i$ is a weak version of the Meta-theorem of the deduction that can be stated as follows: if $\alpha \vdash \beta$ then $\vdash \alpha \rightarrow  \beta$. In the same way, the rule $\vee_e$ can be stated as follows: if $\alpha \vdash \gamma, \beta \vdash \gamma$ then $\alpha \vee \beta \vdash \gamma$. \\


\begin{definition}function of translation (*) of graphs to formulas.

For every proposition $p$ and for every graphic $A$ and $B$
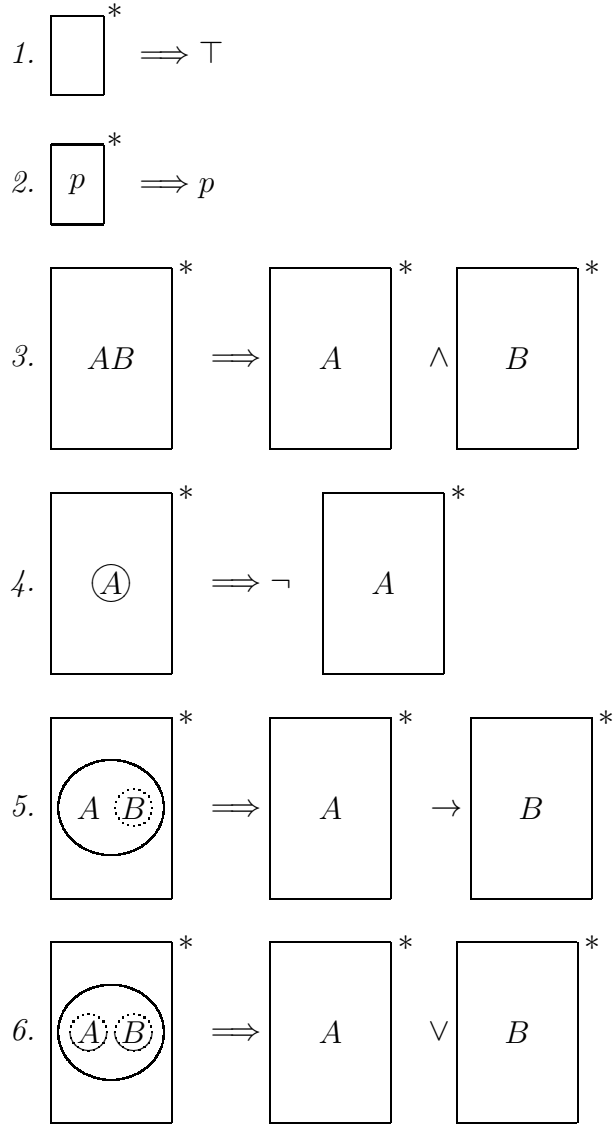
\begin{enumerate}
	\item $\xy(0,0)*=<0.7cm,1.05cm>\frm{-};(5,6)*+{*}\endxy \Longrightarrow \top$
	\item $\xy(0,0)*=<0.7cm,1.05cm>\frm{-}*+{p};(5,6)*+{*}\endxy \Longrightarrow p$
	\item $\xy(0,0)*=<1.6cm,2.4cm>\frm{-}*+{AB};(10,12)*+{*}\endxy \Longrightarrow \xy(0,0)*=<1.6cm,2.4cm>\frm{-}*+{A};(10,12)*+{*}\endxy \wedge \xy(0,0)*=<1.6cm,2.4cm>\frm{-}*+{B};(10,12)*+{*}\endxy$
	\item $\xy(0,0)*=<1.6cm,2.4cm>\frm{-}*+{A}*\cir<7pt>{};(10,12)*+{*}\endxy \Longrightarrow \neg \quad \xy(0,0)*=<1.6cm,2.4cm>\frm{-}*+{A};(10,12)*+{*}\endxy$
	\item $\xy(0,0)*=<1.6cm,2.4cm>\frm{-}*\xycircle<20pt,18pt>{};(3,0)*+{B}*\xycircle<7pt>{.};(-3,0)*+{A};(10,12)*+{*}\endxy \Longrightarrow \xy(0,0)*=<1.6cm,2.4cm>\frm{-}*+{A};(10,12)*+{*}\endxy \rightarrow \xy(0,0)*=<1.6cm,2.4cm>\frm{-}*+{B};(10,12)*+{*}\endxy$
 \item $\xy(0,0)*=<1.6cm,2.4cm>\frm{-}*\xycircle<20pt,18pt>{};(3,0)*+{B}*\xycircle<7pt>{.}*\cir<7pt>{d^u};(-3,0)*+{A}*\xycircle<7pt>{.}*\cir<7pt>{d^u};(10,12)*+{*}\endxy \Longrightarrow \xy(0,0)*=<1.6cm,2.4cm>\frm{-}*+{A};(10,12)*+{*}\endxy \vee \xy(0,0)*=<1.6cm,2.4cm>\frm{-}*+{B};(10,12)*+{*}\endxy$

\end{enumerate}
\end{definition}

\subsubsection{$ALFA_I$ System}

Then, in order to define the $ALFA_I$ system, the rules of the natural deduction system for intuitionism are copied by using the function defined above (as we noted at the beginning of the section), in the same manner the primitive $ALFAo$ symbols are added, it's unique axiom, and the $R_0$ rule.

\newpage
\begin{flushleft}\textbf{Primitive symbols:}\end{flushleft}

\begin{flushleft}
The system's symbols $ALFAo$ and the cuts $\xy*\xycircle<7pt>{.}\endxy$ and $\xy*\cir<7pt>{d^u}*\xycircle<7pt>{.}\endxy$
\end{flushleft}

\begin{definition}
\emph{Well done graphs (wdg):}

\begin{enumerate}
	\item The graphs constructed by the rule of construction of the $ALFAo$ system's graphs are wdgs.
	\item If $
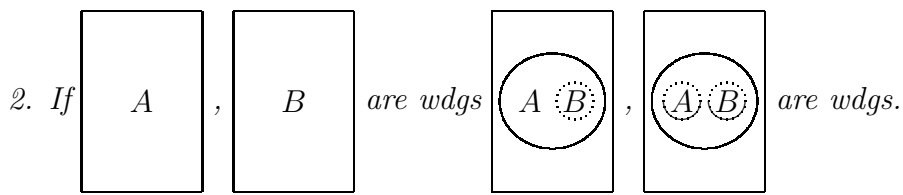
\xy*=<1.6cm,2.4cm>\frm{-}*+{A}\endxy$ , $\xy*=<1.6cm,2.4cm>\frm{-}*+{B}\endxy$ are wdgs $\xy(0,0)*=<1.6cm,2.4cm>\frm{-}*\xycircle<20pt,18pt>{};(3,0)*+{B}*\xycircle<7pt>{.};(-3,0)*+{A}\endxy$ , $\xy(0,0)*=<1.6cm,2.4cm>\frm{-}*\xycircle<20pt,18pt>{};(3,0)*+{B}*\xycircle<7pt>{.}*\cir<7pt>{d^u};(-3,0)*+{A}*\xycircle<7pt>{.}*\cir<7pt>{d^u}\endxy$ are wdgs.
	\item The graphs constructed according to 1. y 2. are wdgs.
\end{enumerate}
\end{definition}


\begin{flushleft}
\textbf{System's axiom}
\end{flushleft}

\xy*=<0.5cm,0.75cm>\frm{-}\endxy\vspace{10pt}

\textbf{System's deductive rules}

\emph{A) First degree deductive rules:}

\begin{flushleft}
Of insertion
\end{flushleft}

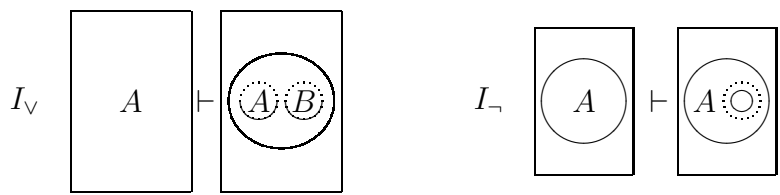
\begin{align*}&
I_{\vee} \quad \xy(-20,0)*=<1.6cm,2.4cm>\frm{-}*+{A};(0,0)*=<1.6cm,2.4cm>\frm{-}*\xycircle<20pt,18pt>{};(3,0)*+{B}*\xycircle<7pt>{.}*\cir<7pt>{d^u};(-3,0)*+{A}*\xycircle<7pt>{.}*\cir<7pt>{d^u};(-10,0)*+{\vdash}\endxy & I_{\neg} \quad \xy(0,0)*=<1.3cm,1.95cm>\frm{-}*\cir<16pt>{}*+{A};(21,0)*\cir<4pt>{}*\xycircle<7pt>{.};(10,0)*+{\vdash};(19,0)*=<1.3cm,1.95cm>\frm{-}*\cir<16pt>{};(16,0)*+{A}\endxy
\end{align*}

\begin{flushleft}
Of elimination (the rule $E_\bot$ can also be considered of insertion)
\end{flushleft}

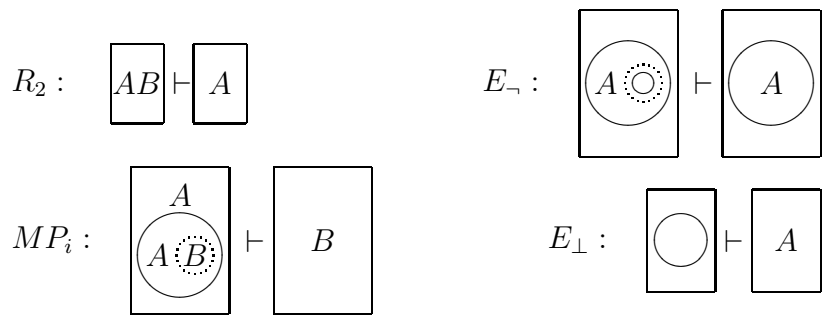
\begin{align*}
& R_2: \quad \xy(0,0)*=<0.7cm,1.05cm>\frm{-}*+{AB};(11,0)*=<0.7cm,1.05cm>\frm{-}*+{A};(6,0)*+{\vdash}\endxy &
E_{\neg}: \quad \xy(0,0)*=<1.3cm,1.95cm>\frm{-}*\cir<16pt>{};(2,0)*\cir<4pt>{}*\xycircle<7pt>{.};(-3,0)*+{A};(10,0)*+{\vdash};(19,0)*=<1.3cm,1.95cm>\frm{-}*+{A}*\cir<16pt>{}\endxy   \\
& MP_i: \quad  \xy(0,0)*=<1.3cm,1.95cm>\frm{-};(2,-2)*+{B}*\xycircle<7pt>{.};(-3,-2)*+{A};(0,6)*+{A};(0,-2)*\cir<16pt>{};(10,0)*+{\vdash};(19,0)*=<1.3cm,1.95cm>\frm{-}*+{B}\endxy
& E_{\bot}: \quad \xy(0,0)*=<0.9cm,1.35cm>\frm{-}*\cir<10pt>{};(14,0)*=<0.9cm,1.35cm>\frm{-}*+{A};(7,0)*+{\vdash}\endxy
\end{align*}

\begin{flushleft}
\emph{B) Second degree deductive rules:}
\end{flushleft}

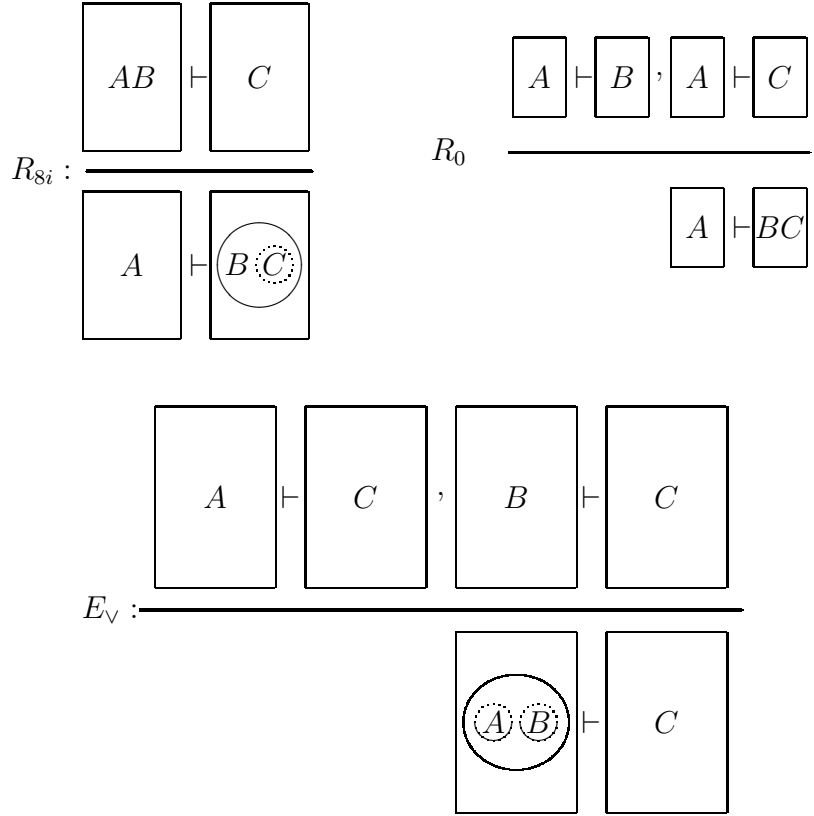
\begin{align*} & \xy(-12,-12.5)*+{R_{8i}:};(0,0)*=<1.3cm,1.95cm>\frm{-}*+{AB};(17,0)*=<1.3cm,1.95cm>\frm{-}*+{C};(9,0)*+{\vdash};(0,-25)*=<1.3cm,1.95cm>\frm{-}*+{A};(17,-25)*=<1.3cm,1.95cm>\frm{-}*\cir<16pt>{};(14,-25)*+{B};(19,-25)*+{C}*\xycircle<7pt>{.};(9,-25)*+{\vdash};(9,-12.5)*=<3cm,0.01cm>\frm{-}\endxy & \xy(-12,-10)*+{R_0};(11,0)*=<0.7cm,1.05cm>\frm{-}*+{B};(0,0)*=<0.7cm,1.05cm>\frm{-}*+{A};(6,0)*+{\vdash};(32,0)*=<0.7cm,1.05cm>\frm{-}*+{C};(21,0)*=<0.7cm,1.05cm>\frm{-}*+{A};(27,0)*+{\vdash};(32,-20)*=<0.7cm,1.05cm>\frm{-}*+{BC};(21,-20)*=<0.7cm,1.05cm>\frm{-}*+{A};(27,-20)*+{\vdash};(16,-10)*=<4cm,0.01cm>\frm{-};(16,0)*+{,}\endxy
\end{align*}

\begin{align*} &
\xy(-14,-15)*+{E_{\vee}:};(20,0)*=<1.6cm,2.4cm>\frm{-}*+{C};(0,0)*=<1.6cm,2.4cm>\frm{-}*+{A};(10,0)*+{\vdash};(60,0)*=<1.6cm,2.4cm>\frm{-}*+{C};(40,0)*=<1.6cm,2.4cm>\frm{-}*+{B};(50,0)*+{\vdash};(30,0)*+{,};(40,-30)*=<1.6cm,2.4cm>\frm{-}*\xycircle<20pt,18pt>{};(43,-30)*+{B}*\xycircle<7pt>{.}*\cir<7pt>{d^u};(37,-30)*+{A}*\xycircle<7pt>{.}*\cir<7pt>{d^u};(60,-30)*=<1.6cm,2.4cm>\frm{-}*+{C};(50,-30)*+{\vdash};(30,-15)*=<8cm,0.01cm>\frm{-}\endxy
\end{align*}


Note: The translation of the $\rightarrow_i$ rule results in a weaker rule than $R_{8i}$. It was decided to take $R_{8i}$ instead of it's weaker similar, to achieve some homogeneity or similarity with $ALFAo$. However, in the appendix, the $R_{8i}$ was tested using the rule that translates directly from $\rightarrow_i$, named as $R_{8id}$.

\subsubsection{Deductions in $ALFA_I$}

\begin{theorem}
$R_{8i}^-:\xy(0,13)*=<1.3cm,1.95cm>\frm{-}*+{A};(20,13)*=<1.3cm,1.95cm>\frm{-}*\cir<16pt>{};(16.5,13)*+{B};(21,13)*+{C}*\xycircle<7.5pt>{.};(11,13)*+{\vdash};(11,0)*=<3.5cm,0.01cm>\frm{-};(0,-13)*=<1.3cm,1.95cm>\frm{-}*+{AB};(20,-13)*=<1.3cm,1.95cm>\frm{-}*+{C};(11,-13)*+{\vdash}\endxy \in ALFA_I$. 
\end{theorem}

The proof of this rule is similar to clasic proof.

\begin{theorem}
$\xy(0,8)*=<0.7cm,1.05cm>\frm{-};(0,8)*+{C};(-11,8)*=<0.7cm,1.05cm>\frm{-};(-11,8)*+{B};(-5,8)*+{\underset{R}{\vdash}};;(-5,0)*=<2cm,0.01cm>\frm{-};;(0,-8)*=<0.7cm,1.05cm>\frm{-};(0,-6)*+{A};(0,-10)*+{C};(-11,-8)*=<0.7cm,1.05cm>\frm{-};(-11,-6)*+{A};(-11,-10)*+{B};(-5,-8)*+{\vdash}\endxy \in ALFA_I$
\end{theorem}

$\xy(16,0)*=<0.7cm,1.05cm>\frm{-};(16,-2)*+{B};(5,0)*=<0.7cm,1.05cm>\frm{-};(5,2)*+{A};(5,-2)*+{B};(11,0)*+{\underset{R_2}{\vdash}};(22,0)*+{\underset{R}{\vdash}};(28,0)*=<0.7cm,1.05cm>\frm{-};(28,-2)*+{C};;(52,0)*=<0.7cm,1.05cm>\frm{-};(52,2)*+{A};(41,0)*=<0.7cm,1.05cm>\frm{-};(41,2)*+{A};(41,-2)*+{B};(47,0)*+{\underset{R_2}{\vdash}};;(28,-8)*=<5.4cm,0.01cm>\frm{-};(-1,-8)*+{R_0};;(52,-16)*=<0.7cm,1.05cm>\frm{-};(52,-14)*+{A};(52,-18)*+{C};(41,-16)*=<0.7cm,1.05cm>\frm{-};(41,-14)*+{A};(41,-18)*+{B};(47,-16)*+{\vdash};(35,0)*+{,}\endxy$\vspace{10pt}

The rule just demonstrated is a direct consequence of the use of $R_0$. In other systems such as $alfa$ and the $GEI$ presented by Arnold Oostra, this rule is assumed in meta-language. Henceforth, the use of it will not be sought in order to simplify the presentation of the demonstrations.

\begin{theorem} $I_{p2}:\xy(0,0)*=<1.6cm,2.4cm>\frm{-}*\cir<20pt>{};(-4,0)*+{B};(2,0)*+{A}*\cir<7pt>{}*\xycircle<10pt>{.};(-11,0)*+{\vdash};(-22,0)*=<1.6cm,2.4cm>\frm{-}*\cir<20pt>{};(-26,0)*+{B};(-20,0)*+{A}\endxy \in ALFA_I$
\end{theorem}

\begin{flushleft}
\xy(0,-25)*=<1.6cm,2.4cm>\frm{-}*\cir<20pt>{};(-4,-25)*+{B};(2,-25)*+{A}*\cir<7pt>{}*\xycircle<10pt>{.};(-11,-25)*+{\vdash};(-22,-25)*=<1.6cm,2.4cm>\frm{-}*\cir<20pt>{};(-26,-25)*+{B};(-20,-25)*+{A};(-22,-10)*=<6cm,0.01cm>\frm{-};(-55,-10)*+{R_{8i}};(0,5)*=<1.6cm,2.4cm>\frm{-}*\cir<20pt>{};(-4,5)*+{A};(-11,5)*+{\underset{E_{\neg}}{\vdash}};(-22,5)*=<1.6cm,2.4cm>\frm{-}*\cir<20pt>{};(-26,5)*+{A};(-20,5)*\cir<7pt>{}*\xycircle<10pt>{.};(-41,20)*=<6cm,0.01cm>\frm{-};(-75,20)*+{R_{8i}};(-44,5)*=<1.6cm,2.4cm>\frm{-}*\cir<20pt>{};(-48,5)*+{B};(-44,5)*+{A};(-34,5)*+{\vdash};(-46,-5)*+{B};(-21,35)*=<1.6cm,2.4cm>\frm{-}*\cir<7pt>{};(-32,35)*+{\underset{MP_i}{\vdash}};(-43,35)*=<1.6cm,2.4cm>\frm{-}*\cir<20pt>{};(-47,35)*+{B};(-44,35)*+{A};(-40,35)*\cir<5pt>{}*\xycircle<7pt>{.};(-47,25)*+{B};(-44,25)*+{A};(-64,35)*=<1.6cm,2.4cm>\frm{-}*\cir<20pt>{};(-68,35)*+{B};(-64,35)*+{A};(-54,35)*+{\underset{I_{\neg}}{\vdash}};(-66,25)*+{B};(-63,25)*+{A}\endxy
\end{flushleft}

\begin{theorem}
$E_p:
\xy(0,0)*=<1.6cm,2.4cm>\frm{-}*\cir<18pt>{};(-3,0)*+{A};(2,0)*+{B}*\cir<7pt>{};(-24,0)*=<1.6cm,2.4cm>\frm{-}*\cir<18pt>{};(-27,0)*+{A};(-22,0)*+{B}*\xycircle<7pt>{.};(-14,0)*+{\vdash}\endxy \in ALFA_I$
\end{theorem}

\begin{flushleft}
\xy(0,5)*=<1.6cm,2.4cm>\frm{-}*\cir<20pt>{};(-4,7)*+{A};(1,7)*+{B}*\cir<7pt>{};(-11,5)*+{\underset{E_{\neg}}{\vdash}};(-22,5)*=<1.6cm,2.4cm>\frm{-}*\cir<20pt>{};(-26,7)*+{A};(-21,7)*+{B}*\cir<7pt>{};(-23,2)*\cir<3pt>{}*\xycircle<5pt>{.};(-41,20)*=<6cm,0.01cm>\frm{-};(-75,20)*+{R_{8i}};(-44,5)*=<1.6cm,2.4cm>\frm{-}*\cir<18pt>{};(-47,5)*+{A};(-42,5)*+{B}*\xycircle<7pt>{.};(-34,5)*+{\vdash};(-21,35)*=<1.6cm,2.4cm>\frm{-};(-21,30)*\cir<5pt>{};(-32,35)*+{\underset{MP_i}{\vdash}};(-43,35)*=<1.6cm,2.4cm>\frm{-};(-44,38)*+{B};(-44,27)*+{B}*\cir<7pt>{};(-66,35)*=<1.6cm,2.4cm>\frm{-};(-66,38)*\cir<18pt>{};(-69,38)*+{A};(-64,38)*+{B}*\xycircle<7pt>{.};(-54,35)*+{\underset{MP_i}{\vdash}};(-68,27)*+{A};(-63,27)*+{B}*\cir<7pt>{}\endxy
\end{flushleft}

\begin{theorem}$I_{p3}: \xy(-10,0)*=<1.6cm,2.4cm>\frm{-}*\xycircle<20pt,18pt>{};(-7,0)*+{B}*\xycircle<7pt>{.}*\cir<7pt>{d^u};(-13,0)*+{A}*\xycircle<7pt>{.}*\cir<7pt>{d^u};(10,0)*=<1.6cm,2.4cm>\frm{-}*\xycircle<20pt,18pt>{};(13,0)*+{B}*\xycircle<7pt>{.};(7,0)*+{A}*\cir<7pt>{};(0,0)*+{\vdash}\endxy \in ALFA_I $
\end{theorem}

\begin{flushleft}
\xy(-22,5)*=<1.6cm,2.4cm>\frm{-}*\cir<20pt>{};(-25,5)*+{A}*\cir<7pt>{};(-19,5)*+{B}*\xycircle<7pt>{.};(-33,20)*=<4.5cm,0.01cm>\frm{-};(-60,20)*+{R_{8i}};(-44,5)*=<1.6cm,2.4cm>\frm{-}*+{B};(-34,5)*+{\vdash};(-21,35)*=<1.6cm,2.4cm>\frm{-}*+{B};(-32,35)*+{\underset{R_2}{\vdash}};(-43,35)*=<1.6cm,2.4cm>\frm{-};(-40,35)*+{B};(-46,35)*+{A}*\cir<7pt>{};(-8,20)*+{,};(48,5)*=<1.6cm,2.4cm>\frm{-}*\cir<20pt>{};(45,5)*+{A}*\cir<7pt>{};(51,5)*+{B}*\xycircle<7pt>{.};(29,20)*=<6cm,0.01cm>\frm{-};(-5,20)*+{R_{8i}};(26,5)*=<1.6cm,2.4cm>\frm{-}*+{A};(36,5)*+{\vdash};(49,35)*=<1.6cm,2.4cm>\frm{-}*+{B};(38,35)*+{\underset{E_{\bot}}{\vdash}};(27,35)*=<1.6cm,2.4cm>\frm{-}*\cir<18pt>{};(4,35)*=<1.6cm,2.4cm>\frm{-};(1,35)*+{A};(6,35)*+{A}*\cir<7pt>{};(16,35)*+{\underset{MP_i}{\vdash}};(0,-10)*=<11cm,0.01cm>\frm{-};(-60,-10)*+{E_{\vee}};(-10,-25)*=<1.6cm,2.4cm>\frm{-}*\xycircle<20pt,18pt>{};(-7,-25)*+{B}*\xycircle<7pt>{.}*\cir<7pt>{d^u};(-13,-25)*+{A}*\xycircle<7pt>{.}*\cir<7pt>{d^u};(10,-25)*=<1.6cm,2.4cm>\frm{-}*\xycircle<20pt,18pt>{};(13,-25)*+{B}*\xycircle<7pt>{.};(7,-25)*+{A}*\cir<7pt>{};(0,-25)*+{\vdash}\endxy
\end{flushleft}

\begin{theorem} $\xy(-10,0)*=<1.6cm,2.4cm>\frm{-}*\xycircle<20pt,18pt>{};(-7,0)*+{B}*\xycircle<7pt>{.}*\cir<7pt>{d^u};(-13,0)*+{A}*\xycircle<7pt>{.}*\cir<7pt>{d^u};(10,0)*=<1.6cm,2.4cm>\frm{-}*\xycircle<20pt,18pt>{};(13,0)*+{B}*\cir<7pt>{};(7,0)*+{A}*\cir<7pt>{};(0,0)*+{\vdash}\endxy \in ALFA_I $
\end{theorem}

$\xy(-10,0)*=<1.6cm,2.4cm>\frm{-}*\xycircle<20pt,18pt>{};(-7,0)*+{B}*\xycircle<7pt>{.}*\cir<7pt>{d^u};(-13,0)*+{A}*\xycircle<7pt>{.}*\cir<7pt>{d^u};(10,0)*=<1.6cm,2.4cm>\frm{-}*\xycircle<20pt,18pt>{};(13,0)*+{B}*\xycircle<7pt>{.};(7,0)*+{A}*\cir<7pt>{};(0,0)*+{\underset{I_{p3}}{\vdash}};(30,0)*=<1.6cm,2.4cm>\frm{-}*\xycircle<20pt,18pt>{};(33,0)*+{B}*\cir<7pt>{};(27,0)*+{A}*\cir<7pt>{};(20,0)*+{\underset{E_{p}}{\vdash}}\endxy$\\

Some of the demonstrated theorems suggest the change of the basic rules of $ALFA_I$ by those, given its simplicity, because it allows to obtain a simpler system to handle and of better geometric appearance. Next we present the $ALFA_{Io}$ system that results from changing the rules $I_\neg, E_\bot, E_\neg$ of $ALFA_I$ for the new theorems $I_{p1}, E_{p2}, I_{p3}$.

\newpage

\subsection{Deductive $ALFA_{Io}$ system}

\textbf{System's axiom}

\xy*=<0.5cm,0.75cm>\frm{-}\endxy\vspace{10pt}

\textbf{System's deductive rules}

\emph{A)Conserved rules:}

\begin{flushleft}
$MP_i, I_\vee, R_{8i}, R_0, E_\vee, R_2$.
\end{flushleft}

\emph{B)New rules:}

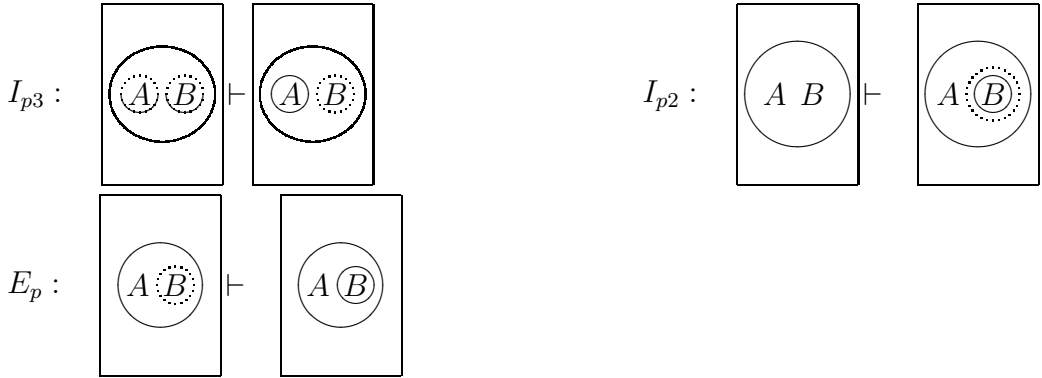
\begin{flalign*}&
I_{p3}: \quad \xy(-10,0)*=<1.6cm,2.4cm>\frm{-}*\xycircle<20pt,18pt>{};(-7,0)*+{B}*\xycircle<7pt>{.}*\cir<7pt>{d^u};(-13,0)*+{A}*\xycircle<7pt>{.}*\cir<7pt>{d^u};(10,0)*=<1.6cm,2.4cm>\frm{-}*\xycircle<20pt,18pt>{};(13,0)*+{B}*\xycircle<7pt>{.};(7,0)*+{A}*\cir<7pt>{};(0,0)*+{\vdash}\endxy &
I_{p2}: \quad
\xy(0,0)*=<1.6cm,2.4cm>\frm{-}*\cir<20pt>{};(-4,0)*+{A};(2,0)*+{B}*\cir<7pt>{}*\xycircle<10pt>{.};(-24,0)*=<1.6cm,2.4cm>\frm{-}*\cir<20pt>{};(-27,0)*+{A};(-22,0)*+{B};(-14,0)*+{\vdash}\endxy \\
& E_p: \quad \xy(0,0)*=<1.6cm,2.4cm>\frm{-}*\cir<18pt>{};(-3,0)*+{A};(2,0)*+{B}*\cir<7pt>{};(-24,0)*=<1.6cm,2.4cm>\frm{-}*\cir<18pt>{};(-27,0)*+{A};(-22,0)*+{B}*\xycircle<7pt>{.};(-14,0)*+{\vdash}\endxy
\end{flalign*}\\


$ALFA_{Io}$ is an $GEI$ system whose rules propose new geomorphic content compared to other existential graphic systems. The $E_p$ rule is a clear example of this, since it can be stated as: \emph{any single dotted closed curve can be closed (complete)}; With which we have to go from an implication of two-graphs to it's equivalent in terms of negation and conjunction, it is enough to close the dotted curves. In an analogous way, it happens with $I_{p3}$, where the passage from the disjunction of two graphs to their implication derives from closing one of the semi-dotted curves and opening (the opposite of closing) the other semi-dotted curve. In this way, a simple management system is obtained with high geomorphic value.

$ALFA_{Io}$ is equivalent to $ALFA_I$, this is proved by demonstrating that $I_\neg, E_\bot, E_\neg$ are deductions of the first.


\subsubsection{Deductions in $ALFA_{Io}$}

\begin{theorem}
$I_c: \quad \xy(0,0)*=<1.3cm,1.95cm>\frm{-}*+{A};(20,0)*=<1.3cm,1.95cm>\frm{-}*\cir<16pt>{}*+{A}*\cir<7pt>{};(10,0)*+{\vdash}\endxy \in ALFA_{Io}$
\end{theorem}

\begin{flushleft}
\xy(40,0)*=<1.3cm,1.95cm>\frm{-}*+{A};(60,0)*=<1.3cm,1.95cm>\frm{-}*+{A};(50,0)*+{\underset{R_2}{\vdash}};(50,-13)*=<3.5cm,0.01cm>\frm{-};(30,-13)*+{R_{8i}};(40,-26)*=<1.3cm,1.95cm>\frm{-}*+{A};(60,-26)*=<1.3cm,1.95cm>\frm{-}*\cir<16pt>{};(61,-26)*+{A}*\xycircle<7pt>{.};(80,-26)*=<1.3cm,1.95cm>\frm{-}*\cir<16pt>{}*+{A}*\cir<7pt>{};(51,-26)*+{\vdash};(70,-26)*+{\underset{E_p}{\vdash}}\endxy
\end{flushleft}

\begin{theorem}$E_\bot \in ALFA_{Io}$
\end{theorem}

\begin{flushleft} \xy(-20,0)*=<1.6cm,2.4cm>\frm{-}*\cir<9pt>{};(0,0)*=<1.6cm,2.4cm>\frm{-};(0,0)*\xycircle<20pt,18pt>{};(3,0)*+{A}*\xycircle<7pt>{.}*\cir<7pt>{d^u};(-3,0)*\cir<5pt>{}*\xycircle<7pt>{.}*\cir<7pt>{d^u};(-10,0)*+{\underset{I_{\vee}}{\vdash}};;(20,0)*=<1.6cm,2.4cm>\frm{-};(20,0)*\xycircle<20pt,18pt>{};(23,0)*+{A}*\xycircle<7pt>{.};(17,0)*\cir<5pt>{}*\cir<7pt>{};(10,0)*+{\underset{I_{p3}}{\vdash}};;(40,0)*=<1.6cm,2.4cm>\frm{-};(40,0)*\xycircle<20pt,18pt>{};(43,0)*+{A}*\xycircle<7pt>{.};(37,0)*\cir<7pt>{}*\cir<5pt>{};(45,9)*\cir<7pt>{}*\cir<5pt>{};(30,0)*+{\underset{I_c}{\vdash}};;(62,0)*=<1.6cm,2.4cm>\frm{-}*+{A};(51,0)*+{\underset{MP_i}{\vdash}}\endxy
\end{flushleft}

\begin{theorem}$R_5^{'} \in ALFA_{Io}$ 
\end{theorem}

$R_5^{'}$ is a particular case of $R_5$, rule that is spoken of in the first section.

\xy(-22,5)*=<1.6cm,2.4cm>\frm{-}*\cir<20pt>{};(-26,5)*+{A};(-22,15)*+{A};(-20,5)*+{B}*\cir<7pt>{}*\xycircle<10pt>{.};(-44,5)*=<1.6cm,2.4cm>\frm{-}*\cir<20pt>{}*+{AB};(-44,15)*{A};(-34,5)*+{\underset{I_{p2}}{\vdash}};(0,5)*=<1.6cm,2.4cm>\frm{-}*+{B}*\cir<7pt>{};(-11,5)*+{\underset{MP_i}{\vdash}}\endxy

\begin{theorem}$I_\neg \in ALFA_{Io}$
\end{theorem}

\xy(0,0)*=<1.3cm,1.95cm>\frm{-}*\cir<16pt>{}*+{A};(10,0)*+{\underset{I_{p2}}{\vdash}};(19,0)*=<1.3cm,1.95cm>\frm{-}*\cir<16pt>{};(17,0)*+{A};(21,0)*\cir<5pt>{}*\xycircle<7pt>{.}\endxy

\begin{theorem}$E_\neg \in ALFA_{Io}$
\end{theorem}

\begin{flushleft} \xy(0,0)*=<1.3cm,1.95cm>\frm{-};(2,-2)*\cir<5pt>{}*\xycircle<7pt>{.};(-3,-2)*+{A};(0,-2)*\cir<16pt>{};(10,0)*+{\underset{E_p}{\vdash}};(19,0)*=<1.3cm,1.95cm>\frm{-};(21,-2)*\cir<5pt>{}*\cir<7pt>{};(16,-2)*+{A};(19,-2)*\cir<16pt>{};;(29,0)*+{\underset{I_c}{\vdash}};(38,0)*=<1.3cm,1.95cm>\frm{-};(40,-2)*\cir<5pt>{}*\cir<7pt>{};(35,-2)*+{A};(38,-2)*\cir<16pt>{};(41,6.5)*\cir<7pt>{}*\cir<5pt>{};;(48,0)*+{\underset{R_5^{'}}{\vdash}};(57,0)*=<1.3cm,1.95cm>\frm{-};(54,-2)*+{A};(57,-2)*\cir<16pt>{}\endxy
\end{flushleft}

The previous theorems prove the equivalence between $ALFA_{Io}$ and $ALFA_I$. Now it remains to answer a question: Is it possible to obtain a system of existential graphs equivalent to the $CPC$ of which $ALFA_{Io}$ is a sub-system of?

We will then deduce those $ALFAo$ rules that belong to $ALFA_{Io}$.

\begin{theorem}$R_5 \in ALFA_{Io}$
\end{theorem}

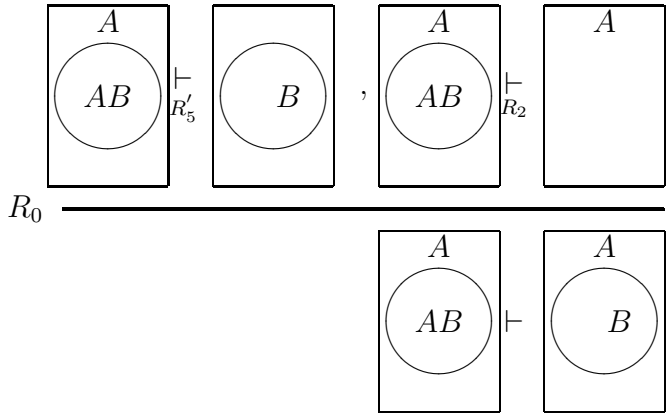
\begin{flushleft}
\xy(-22,5)*=<1.6cm,2.4cm>\frm{-}*\cir<20pt>{};(-20,5)*+{B};(-44,5)*=<1.6cm,2.4cm>\frm{-}*\cir<20pt>{}*+{AB};(-44,15)*{A};(-34,5)*+{\underset{R_5^{'}}{\vdash}};;(22,5)*=<1.6cm,2.4cm>\frm{-};(22,15)*+{A};(0,5)*=<1.6cm,2.4cm>\frm{-}*\cir<20pt>{}*+{AB};(0,15)*{A};(10,5)*+{\underset{R_2}{\vdash}};;(22,-25)*=<1.6cm,2.4cm>\frm{-}*\cir<20pt>{};(22,-15)*+{A};(24,-25)*+{B};(0,-25)*=<1.6cm,2.4cm>\frm{-}*\cir<20pt>{}*+{AB};(0,-15)*{A};(10,-25)*+{\vdash};(-55,-10)*+{R_0};(-10,-10)*=<8cm,0cm>\frm{-};(-10,5)*+{,}\endxy
\end{flushleft}

\begin{theorem}$R_7 \in ALFA_{Io}$
\end{theorem}

$
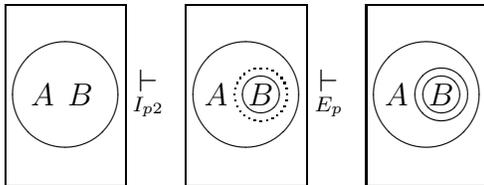
\xy(0,0)*=<1.6cm,2.4cm>\frm{-}*\cir<20pt>{};(-4,0)*+{A};(2,0)*+{B}*\cir<7pt>{}*\xycircle<10pt>{.};(-24,0)*=<1.6cm,2.4cm>\frm{-}*\cir<20pt>{};(-27,0)*+{A};(-22,0)*+{B};(-13,0)*+{\underset{I_{p2}}{\vdash}};(24,0)*=<1.6cm,2.4cm>\frm{-}*\cir<20pt>{};(20,0)*+{A};(26,0)*+{B}*\cir<7pt>{}*\cir<10pt>{};(11,0)*+{\underset{E_p}{\vdash}}\endxy$

$R_7^-$ is deduced in order to facilitate the deduction of the $R_3$ rule.

\begin{theorem}
$R_7^- \in ALFA_{Io}$
\end{theorem}

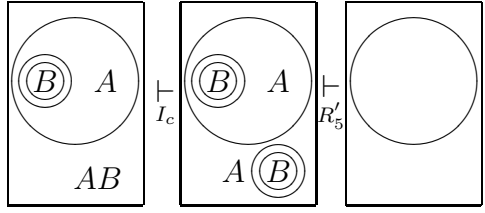
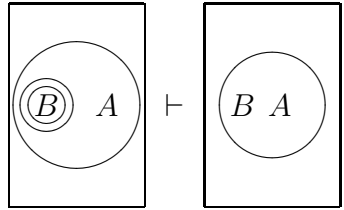
\xy(-22,5)*=<1.8cm,2.7cm>\frm{-}*\cir<23pt>{};(-26,5)*+{B};(-21,5)*+{A};(-41,25)*=<6cm,0.01cm>\frm{-};(-80,25)*+{R_{8i}};(-48,5)*=<1.8cm,2.7cm>\frm{-}*\cir<24pt>{};(-52,5)*+{B}*\cir<7pt>{}*\cir<10pt>{};(-44,5)*+{A};(-35,5)*+{\vdash};(-21,45)*=<1.8cm,2.7cm>\frm{-};(-21,48)*\cir<24pt>{};(-32,45)*+{\underset{R_5^{'}}{\vdash}};(-43,45)*=<1.8cm,2.7cm>\frm{-};(-43,48)*\cir<24pt>{};(-47,48)*+{B}*\cir<7pt>{}*\cir<10pt>{};(-39,48)*+{A};(-39,36)*+{B}*\cir<7pt>{}*\cir<10pt>{};(-45,36,)*+{A};(-66,45)*=<1.8cm,2.7cm>\frm{-};(-66,48)*\cir<24pt>{};(-70,48)*+{B}*\cir<7pt>{}*\cir<10pt>{};(-62,48)*+{A};(-54,45)*+{\underset{I_c}{\vdash}};(-63,35)*+{AB}\endxy

\begin{theorem}$R_3 \in ALFA_{Io}$
\end{theorem}

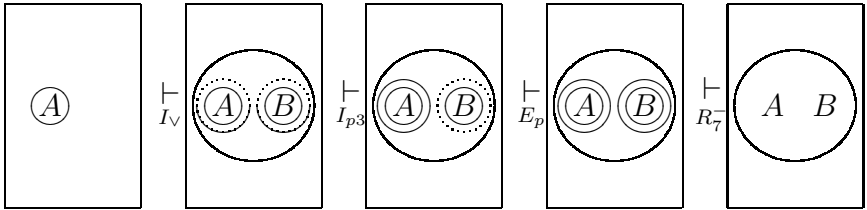
\xy(-30,0)*=<1.8cm,2.7cm>\frm{-};(-33,0)*+{A}*\cir<7pt>{};(-17,0)*+{\underset{I_\vee}{\vdash}};;(-6,0)*=<1.8cm,2.7cm>\frm{-}*\xycircle<23pt,21pt>{};(-2,0)*+{B}*\xycircle<10pt>{.}*\cir<10pt>{d^u}*\cir<7pt>{};(-10,0)*+{A}*\xycircle<10pt>{.}*\cir<10pt>{d^u}*\cir<7pt>{};;(18,0)*=<1.8cm,2.7cm>\frm{-}*\xycircle<23pt,21pt>{};(22,0)*+{B}*\xycircle<10pt>{.}*\cir<7pt>{};(14,0)*+{A}*\cir<7pt>{}*\cir<10pt>{};(7,0)*+{\underset{I_{p3}}{\vdash}};;(42,0)*=<1.8cm,2.7cm>\frm{-}*\xycircle<23pt,21pt>{};(46,0)*+{B}*\cir<10pt>{}*\cir<7pt>{};(38,0)*+{A}*\cir<7pt>{}*\cir<10pt>{};(31,0)*+{\underset{E_p}{\vdash}};;(66,0)*=<1.8cm,2.7cm>\frm{-}*\xycircle<23pt,21pt>{};(70,0)*+{B};(63,0)*+{A};(55,0)*+{\underset{R_7^-}{\vdash}}\endxy\vspace{10pt}

In summary, the rules $R_0,R_2,R_3,R_5$ and $R_7$ belong to both systems ($ALFAo$ and $ALFA_{Io}$).


\subsection{$ALFA_{Io}$ and the $CPC$}

Most axiom systems that formalize the $IPC$, are sub-systems of a system that formalizes the $CPC$\footnote{See examples of it in \cite{Chagrov}}, fact that can be of great utility when one wishes to compare the intuitionist and classical logic. Based on this, we intend to find a system of existential graphs equivalent to $ALFAo$ of which $ALFA_{Io}$ is sub-system; For this it is necessary to add some rules to $ALFA_{Io}$ in such a way that the new system deducts all of $ALFAo$ rules.\\

In order to obtain a classic system from the $ALFA_{Io}$ system, the $I_{\vee p}: \quad\\ 
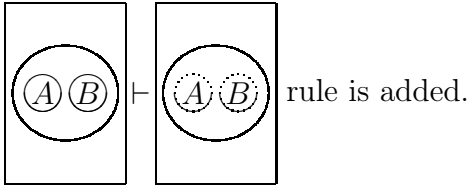
\xy(-10,0)*=<1.6cm,2.4cm>\frm{-}*\xycircle<20pt,18pt>{};(-7,0)*+{B}*\cir<7pt>{}*\cir<7pt>{d^u};(-13,0)*+{A}*\cir<7pt>{}*\cir<7pt>{d^u};(10,0)*=<1.6cm,2.4cm>\frm{-}*\xycircle<20pt,18pt>{};(13,0)*+{B}*\xycircle<7pt>{.}*\cir<7pt>{d^u};(7,0)*+{A}*\xycircle<7pt>{.}*\cir<7pt>{d^u};(0,0)*+{\vdash}\endxy $ rule is added.\\

Now we must prove that the following equivalences are in the $ALFA_{Io}+\{I_{\vee p}\}$ system:\\

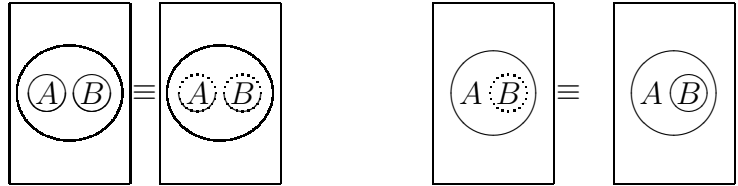
\begin{align*}&
\xy(-10,0)*=<1.6cm,2.4cm>\frm{-}*\xycircle<20pt,18pt>{};(-7,0)*+{B}*\cir<7pt>{}*\cir<7pt>{d^u};(-13,0)*+{A}*\cir<7pt>{}*\cir<7pt>{d^u};(10,0)*=<1.6cm,2.4cm>\frm{-}*\xycircle<20pt,18pt>{};(13,0)*+{B}*\xycircle<7pt>{.}*\cir<7pt>{d^u};(7,0)*+{A}*\xycircle<7pt>{.}*\cir<7pt>{d^u};(0,0)*+{\equiv}\endxy  &\xy(0,0)*=<1.6cm,2.4cm>\frm{-}*\cir<18pt>{};(-3,0)*+{A};(2,0)*+{B}*\cir<7pt>{};(-24,0)*=<1.6cm,2.4cm>\frm{-}*\cir<18pt>{};(-27,0)*+{A};(-22,0)*+{B}*\xycircle<7pt>{.};(-14,0)*+{\equiv}\endxy
\end{align*}\\

The first equivalence is given so $I_{p3} \in ALFA_{Io}$ and $I_{\vee p}$ is the added rule. To prove the second equivalence, it is enough to prove that the $E_p^{-1}$ rule belongs to the new system, since $E_p \in ALFA_{Io}$; which is shown below:

\begin{theorem}$E_p^{-1} \in ALFA_{Io}+\{I_{\vee p}\}$
\end{theorem}

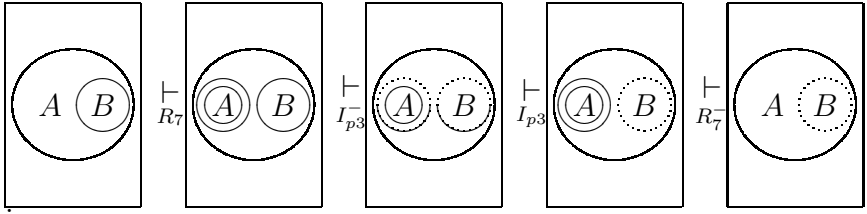
\xy(-30,0)*=<1.8cm,2.7cm>\frm{-}*\xycircle<23pt,21pt>{};(-33,0)*+{A};(-26,0)*+{B}*\cir<10pt>{};(-17,0)*+{\underset{R_7}{\vdash}};;(-6,0)*=<1.8cm,2.7cm>\frm{-}*\xycircle<23pt,21pt>{};(-2,0)*+{B}*\cir<10pt>{};(-10,0)*+{A}*\cir<10pt>{}*\cir<7pt>{};;(18,0)*=<1.8cm,2.7cm>\frm{-}*\xycircle<23pt,21pt>{};(22,0)*+{B}*\xycircle<10pt>{.}*\cir<10pt>{d^u};(14,0)*+{A}*\xycircle<10pt>{.}*\cir<10pt>{d^u}*\cir<7pt>{};(7,0)*+{\underset{I_{p3}^-}{\vdash}};;(42,0)*=<1.8cm,2.7cm>\frm{-}*\xycircle<23pt,21pt>{};(46,0)*+{B}*\xycircle<10pt>{.};(38,0)*+{A}*\cir<7pt>{}*\cir<10pt>{};(31,0)*+{\underset{I_{p3}}{\vdash}};;(66,0)*=<1.8cm,2.7cm>\frm{-}*\xycircle<23pt,21pt>{};(70,0)*+{B}*\xycircle<10pt>{.};(63,0)*+{A};(55,0)*+{\underset{R_7^-}{\vdash}}\endxy.\\

After seeing that in the new system these equivalences are maintained, it's natural to ask if it's really classic. To see that it is, we will prove that $ALFA_{Io}+\{I_{\vee p}\} \equiv ALFAo$, for which it should only be shown that all rules of $ALFAo$ are theorems of $ALFA_{Io}+\{I_{\vee p}\}$.

\begin{theorem}$R_6 \in ALFA_{Io}+\{I_{\vee p}\}$
\end{theorem}

\xy(-22,0)*=<1.6cm,2.4cm>\frm{-}*\cir<20pt>{}*\xycircle<15pt>{.}*+{A};(-44,0)*=<1.6cm,2.4cm>\frm{-}*\cir<20pt>{}*\cir<15pt>{}*+{A};(-34,0)*+{\underset{E_p^-}{\vdash}};(0,0)*=<1.6cm,2.4cm>\frm{-}*+{B};(-11,0)*+{\underset{MP_i}{\vdash}}\endxy

In the following deduction  $R_6$ is used as theorem of the $ALFA_{Io}+\{I_{\vee p}\}$ system.

\begin{theorem}$R_4 \in ALFA_{Io}+\{I_{\vee p}\}$
\end{theorem}

\begin{flushleft}
\xy(24,5)*=<1.6cm,2.4cm>\frm{-}*\cir<20pt>{};(20,8)*+{A};(20,4)*+{C};(26,5)*+{BC}*\cir<11pt>{};(12,5)*+{\underset{E_p}{\vdash}};;(1,5)*=<1.6cm,2.4cm>\frm{-}*\cir<20pt>{};(-3,8)*+{A};(-3,4)*+{C};(3,5)*+{BC}*\xycircle<11pt>{.};(-32,20)*=<8cm,0.01cm>\frm{-};(-74,20)*+{R_{8i}};(-21,5)*=<1.6cm,2.4cm>\frm{-}*\cir<20pt>{};(-25,7)*+{A};(-25,3)*+{C};(-19,5)*+{B}*\cir<7pt>{};(-11,5)*+{\vdash};;(1,35)*=<1.6cm,2.4cm>\frm{-};(3,38)*+{B};(4,27)*+{C};(-10,35)*+{\underset{R_2}{\vdash}};;(-21,35)*=<1.6cm,2.4cm>\frm{-};(-19,38)*+{B};(-18,27)*+{C};(-23,27)*+{A};;(-32,35)*+{\underset{R_6}{\vdash}};(-43,35)*=<1.6cm,2.4cm>\frm{-};(-40,27)*+{C};(-45,27)*+{A};(-43,38)*\cir<20pt>{};(-41,38)*+{B}*\cir<7pt>{};;(-66,35)*=<1.6cm,2.4cm>\frm{-};(-66,38)*\cir<20pt>{};(-70,39)*+{A};(-70,35)*+{C};(-64,38)*+{B}*\cir<7pt>{};(-54,35)*+{\underset{R_5}{\vdash}};(-68,27)*+{A};(-63,27)*+{C}\endxy
\end{flushleft}

The rules $R_0,R_2,R_3,R_5$ y $R_7$ belong to $ALFA_{Io}$, therefore, also to $ALFA_{Io}+\{I_{\vee p}\}$, and the rules $R_4$ and $R_6$ belong to $ALFA_{Io}+\{I_{\vee p}\}$. It remains to be seen that $R_8$ is deductible from this system, proving that it is immediate considering the equivalences between the graphs developed previously.

In conclusion, $ALFA_{Io}+\{I_{\vee p}\}$ and $ALFAo$ are equivalent, except that the first one has more symbols than the second.

\section{The $GEI$ of Arnold Oostra}

Arnold Oostra presents a system of intuitionist existential graphs in \cite{Oostra}. In this article Oostra introduces the system of graphs mentioned as a proposal for formalization of the intuitionist logic, using diagrams that appear in the manuscripts of Charles S. Peirce. This $GEI$ system differs from the one presented in the previous chapter to a large extent, as will be seen in the following section. It is recommended to refer to the article to understand in depth the proposal made by Oostra.

\subsection{Differences between $ALFA_{Io}$ and the $GEI-Oostra$}

Earlier it had been suggested that $ALFA_{Io}$ could be considered as the intuitionist version of $ALFAo$ and $GEI-Oostra$ as the intuitionistic version of $alfa$. Thus, it is expected that the differences between $ALFAo$ and $alfa$ will also be preserved as differences of $ALFA_{Io}$ with respect to $GEI-Oostra$. In fact, some are maintained, others are not: for the first case, no deductive rule of $ALFA_{Io}$ is defined in function of the parity of the cuts, and the symmetries of the $ALFAo_I$ rules is understood in a different way then the $GEI-Oostra$ system; In the second, $ALFAo_I$ rules are not sub-rules of the $GEI-Oostra$ rules.\\

In addition to the previous differences, there are the differences of the primitive symbols of each system:\\

The $ALFAo_I$ primitive symbols are the ones of $ALFAo$ and the cuts $\xy*\xycircle<7pt>{.}\endxy$ y $\xy*\cir<7pt>{d^u}*\xycircle<7pt>{.}\endxy$\\

The $GEI-Oostra$ primitive symbols are the ones of the $alfa$ system, the curls $\xy(0,0)*\xycircle<30pt,20pt>{};(5,0)*\xycircle<15pt,10pt>{}\endxy$ and the loops $\xy(0,0)*\xycircle<37.5pt,25pt>{};(9.5,0)*\cir<10pt>{};(-9.5,0)*\cir<10pt>{};(0,5)*\cir<10pt>{}\endxy$ \\

Derived from these choices of the primitive symbols of each system, come the differences corresponding to the representations for implication and disjunction:\\

Implication and disjunction of two graphs in $ALFA_{Io}$: $\xy(0,0)*\xycircle<20pt,18pt>{};(3,0)*\xycircle<7pt>{.}\endxy$ , $\xy*\xycircle<20pt,18pt>{};(3,0)*\xycircle<7pt>{.}*\cir<7pt>{d^u};(-3,0)*\xycircle<7pt>{.}*\cir<7pt>{d^u}\endxy$.\\

Implication and disjunction of two graphs in $GEI-Oostra$: $\xy(0,0)*\xycircle<30pt,20pt>{};(-5,0);(5,0)*\xycircle<15pt,10pt>{}\endxy$ , $\xy(0,0)*\xycircle<37.5pt,25pt>{};(9.5,0)*\cir<10pt>{};(-9.5,0)*\cir<10pt>{}\endxy$.\\

$R_o$ is a basic $ALFA_{Io}$ rule stated explicitly as the system's rule, meanwhile in $GEI-Oostra$ no, because in this,the rule is used in the meta-lenguage.\\

In conclusion, an $ALFAo$ style $GEI$ system was obtained that satisfies all expectations desired in this text. Now all that remains is to continue the work by investigating: the intermediate logics and the existential graphs; the extension of $ALFAo$ and $ALFA_{Io}$ to the predicate calculus and the modal logic in the sense that $beta$ and $gamma$ are $alfa$; and to study $ALFA_{Io}$ in a topological sense taking into account the relations of the intuitionist logic and the topology, considering the cuts dotted as open and the rules interpreted as the calculus of the closure of a cut. In this context, it is expected, with the conclusion of this work, to generate more questions and problems that were attempted to solve.

\section{Acknowledgments}
We thank the office of the Vice President of Research at the Technological University of Pereira for funding this research through the project, 3-17-2.

\end{document}